\documentclass{article}

\usepackage[preprint]{nips_2018}

\usepackage[utf8]{inputenc} 
\usepackage[T1]{fontenc}    
\usepackage{hyperref}       
\usepackage{url}            
\usepackage{booktabs}       
\usepackage{amsfonts}       
\usepackage{nicefrac}       
\usepackage{microtype}      

\usepackage{amsthm,amsmath,amssymb}
\usepackage{graphicx}
\usepackage{tabularx}
\usepackage{multirow}
\usepackage{amssymb}
\usepackage{tikz}
\usetikzlibrary{shapes,arrows,positioning}

\usepackage{paralist}
\RequirePackage{hyperref}

\newtheorem{theorem}{Theorem}

\newtheorem{lemma}[theorem]{Lemma} 

\newtheorem{corollary}[theorem]{Corollary}

\theoremstyle{definition} 

\theoremstyle{remark} 

\newtheorem*{remark}{Remark}

\numberwithin{equation}{section}
\numberwithin{theorem}{section}
\numberwithin{example}{section}
\numberwithin{definition}{section}

\numberwithin{figure}{section}

\DeclareMathOperator{\diag}{diag}

\newcommand{\figref}[1]{Figure~\ref{fig:#1}}

\newcommand{\secref}[1]{Section~\ref{sec:#1}}

\newcommand{\lemref}[1]{Lemma~\ref{lem:#1}}

\newcommand{\corref}[1]{Corollary~\ref{cor:#1}}

\title{Algebraic tests of general Gaussian latent tree models}

\author{
  Dennis Leung\\
  Department of Data Sciences and Operations\\
  University of Southern California\\
  \texttt{dmhleung@uw.edu} \\
  \And
  Mathias Drton\\
  Department of Statistics,  University of Washington \&\\
 Department of Mathematical Sciences, University of Copenhagen\\
 \texttt{md5@uw.edu} \\
}

\begin{document}

\maketitle

\begin{abstract}
We consider general Gaussian latent tree models in which the
  observed variables are not restricted to be leaves of the tree.
  Extending related recent work, we give a full semi-algebraic
  description of the set of covariance matrices of any such model.  In
  other words, we find polynomial constraints that characterize when
  a matrix is the covariance matrix of a distribution in a given latent
  tree model.  However, leveraging these constraints to test a given such model is often complicated by
  the number of constraints being large and by singularities of
  individual polynomials, which may invalidate standard approximations
  to relevant probability distributions. Illustrating with the star tree, we propose a new testing methodology
that circumvents singularity issues by
  trading off some statistical estimation efficiency and handles cases with many
  constraints through recent advances on Gaussian approximation for
  maxima of sums of high-dimensional random vectors.
  Our test avoids the need to maximize the possibly multimodal
  likelihood function of such models and is applicable to models with larger number of
  variables.  These points are illustrated in numerical experiments.
\end{abstract}

\section{Introduction}

Latent tree models are associated to a tree-structured graph in which
some nodes represent observed variables and others represent
unobserved (latent) variables.  Due to their tractability, these
models have found many applications in fields ranging from the
traditional life sciences, biology and psychology to contemporary
areas such as artificial intelligence and computer vision; refer to \citet{surveyLT} for a comprehensive review. In this paper, we study the problem of testing the
goodness-of-fit of a postulated Gaussian latent tree model to an
observed dataset.  In a low dimensional setting where the number of
observed variables is small relative to the sample size at hand,
testing is usually based on the likelihood ratio which measures the
divergence in maximum likelihood between the postulated latent tree
model and an unconstrained Gaussian model.  This, however, requires
maximization of the possibly multimodal likelihood function of latent
tree models.  In contrast, recent work of \citet{Shiers} takes a
different approach and leverages known polynomial constraints on the
covariance matrix of the observed variables in a given Gaussian latent
tree.  Specifically, the postulated latent tree is tested with an
aggregate statistic formed from estimates of the polynomial
quantities involved.  This approach can be traced back to
\citet{10.2307/1412107} and \cite{wishart1928sampling}; also see
\citet{algfac, MorMinors}.

We make the following new contributions.  In \secref{char}, we extend
the polynomial characterization of \citet{Shiers} to cases where
observed nodes may also be inner nodes of the tree as considered, for
example, in the tree learning algorithms of \citet{Willsky2013}.
\secref{test} describes how we may use polynomial equality constraints
to test a star tree model.  We base
ourselves on the recent groundbreaking work of \citet{CCKHD}, form our
test statistic as the maximum of unbiased estimates of the relevant
polynomials, and calibrate the critical value for testing based on
multiplier bootstrapping techniques.  This new way of using the polynomials to furnish a test  allows us to handle latent trees with a larger number of observed variables and avoids potential singularity issues caused by individual polynomials.  Numerical experiments in \secref{num-exp} makes comparisons to the likelihood ratio
test and assesses the size of our tests in finite samples.  \secref{conclude} discusses future research directions. 
 \paragraph*{Notation.}
 Let $1\le r\le m$ be two positive integers.  We let
 $[m]=\{1, \dots, m\}$ and
 write ${m \brace r} := \{I \subseteq [m]: |I| = r\}$ for the collection
 of subsets of $[m]$ with cardinality $r$. The supremum norm of a
 vector is written $\|\cdot\|_\infty$. For two random variables $R_1$
 and $R_2$, the symbols $R_1 =_d R_2$ indicates that $R_1$ and $R_2$
 have the same distributions, and $R_1 \approx_d R_2$ indicates that
 the distributions are approximately equal. $N(\mu, \sigma^2)$ means a normal distribution with mean $\mu$ and standard deviation $\sigma^2$. 


%

\section{Characterization of general Gaussian latent trees} \label{sec:char}

We first provide the definition of the models considered in this
paper.  A tree is an undirected graph in which any two nodes are
connected by precisely one path.  Let $T = (V, E)$ be a tree, where
$V$ is the set of nodes, and $E$ is the set of edges which we take to
be unordered duples of nodes in $V$.  We say that $T$ is a latent tree
if it is paired with a set ${\bf X} = \{X_1, \dots, X_m\}\subset V$,
corresponding to $m$ observed variables, such that $v \in {\bf X}$
whenever $v \in V$ is of degree less than or equal to two.  In
particular, $\bf X$ contains all leaf nodes of the tree $T$ (i.e.,
nodes of degree 1), but it may contain additional nodes.  The nodes in
$V \setminus {\bf X}$ correspond to latent variables that are not
observed but each have at least three other neighbors in the tree.
This minimal degree requirements of $3$ on the latent nodes ensures
identifiability \citep[p.1778]{Willsky2013}.  In the terminology of
mathematical phylogenetics, $T$ is a semi-labeled tree on ${\bf X}$
with an injective labeling map; see \citet[p.16]{phylo}.
However, \emph{phylogenetic trees} are latent trees restricted to have
$\bf X$ equal to the set of leaves.  While we have defined $\bf X$ as
a set of nodes, it will be convenient to abuse notation slightly and
let $\bf X$ also denote a random vector $(X_1, \dots, X_m)'$ whose
coordinates correspond to the nodes in question.  The context will
clarify whether we refer to nodes or random variables.

Now we present the polynomial characterization of a Gaussian latent
tree graphical model that extends the results in \citet{Shiers}.  The
Gaussian graphical model on $T$, denoted $\mathcal{M}(T)$, is the set
of all $|V|$-variate Gaussian distributions respecting the pairwise
Markov property of $T$, i.e., for any pair $u, v \in V $ with
$(u, v) \not \in E$, the random variables associated to $u$ and $v$
are conditionally independent given the variables corresponding to
$V \backslash \{u, v\}$.  The $T$-Gaussian latent tree model on
${\bf X}$, denoted $\mathcal{M}_{\bf X}(T)$, is the set of all
$m$-variate Gaussian distributions that are the marginal distribution
for ${\bf X}$ under some distribution in $\mathcal{M}(T)$.  For a
given distribution in $\mathcal{M}(T)$, let $\rho_{pq}$ be the Pearson
correlation of the pair $(X_p, X_q)$ for any
$1 \leq p \not = q \leq m$.  The pairwise Markov property implies that
\begin{equation} \label{para}
 \rho_{pq} = \prod_{(u, v) \in ph_T(X_p, X_q)} \rho'_{uv},
\end{equation}
where $ph_T(X_p, X_q)$ denotes the set of edges on the unique path
that connects $X_p$ and $X_q$ in $T$, and $\rho_{uv}'$ is the Pearson
correlation between a pair of nodes $u$ and $v$ in $V$.  Of course,
$\rho_{uv}'=\rho_{pq}$ if $u=X_p$ and $v=X_q$.  In the sequel, we
often abbreviate $ph_T(X_p, X_q)$ as $ph_T(p, q)$ for simplicity.

Suppose $\Sigma = (\sigma_{pq})_{1 \leq p , q \leq m}$ is the
covariance matrix of ${\bf X}$. Our task is to test whether $\Sigma$
comes from $\mathcal{M}_{\bf X}(T)$ against a saturated Gaussian
graphical model.  We assume that all edges in the tree $T$ correspond
to a nonzero correlation, so that $\Sigma$ contains no zero entries.
The covariance matrices for $\mathcal{M}_{\bf X}(T)$ are parametrized
via~\eqref{para}.  As shown in \cite{Shiers}, this set of covariance
matrices may be characterized by leveraging results on pseudo-metrics
defined on ${\bf X}$.  Suppose
$w : E \longrightarrow \mathbb{R}_{\geq 0}$ is a function that assigns
non-negative weights to the edges in $E$. One can then define a
pseudo-metric
$\delta_w : {\bf X} \times {\bf X} \longrightarrow \mathbb{R}_{\geq
  0}$ by
\[
\delta_w(X_p, X_q) = \left\{
  \begin{array}{lr}
    \sum_{e \in ph_T(p,  q)} w(e) & : p \not= q,\\
    0 & :  p = q.
  \end{array}
\right. 
 \]
 This is known as a $T$-induced pseudo-metric on ${\bf X}$.  The
 following lemma characterizes all the pseudo-metrics on ${\bf X}$
 that are $T$-induced. The proof is a bit delicate and is given in our supplementary material.
\begin{lemma} \label{lem:pseudo}
Suppose $\delta : {\bf X} \times {\bf X} \longrightarrow \mathbb{R}_{\geq 0}$ is a pseudo-metric defined on $\bf X$. Let $\delta_{pq} = \delta(X_p, X_q)$ for any $p, q \in [m]$ for simplicity. Then $\delta$ is a $T$-induced pseudo-metric if and only if   for any four distinct $1 \leq p, q, r, s \leq m$ such that ${ph}_T(p, q)\cap ph_T(r, s) = \emptyset$, 
\begin{equation} \label{T4PC}
\delta_{pq} + \delta_{rs } \leq \delta_{pr} + \delta_{qs}  = \delta_{ps} + \delta_{qr}, 
\end{equation}
and for any three distinct $1 \leq p, q, r \leq m$, 
\begin{equation} \label{T3PC}
\delta_{pq} +  \delta_{qr} = \delta_{pr}
\end{equation}
if $ph_T(p, r) = ph_T(p, q) \cup ph_T(q, r)$. 
\end{lemma}

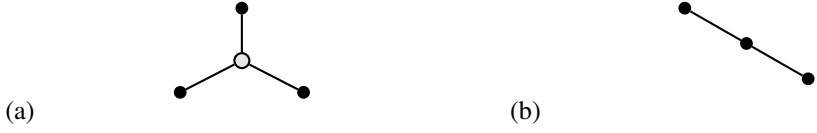
\begin{figure}[t]
  \begin{tabular}{ll}
\begin{minipage}[t]{0.45\hsize}
\centering
  \begin{tikzpicture}[baseline=(1), ->,>=triangle 45, shorten >=0pt,
        auto,thick,
        obsnode/.style={circle,inner sep=0.05cm,fill=black,draw,font=\bf, text= white}, 
 hidnode/.style={circle,inner sep=2pt,fill=gray!20,draw,font=\sffamily},
    label/.style={font=\bf}]
        \node[obsnode,rounded corners] (1) {};

\node[obsnode,rounded corners] [below left = 1cm and 0.7cm of 1] (2) {};

\node[obsnode,rounded corners] [below right = 1cm and 0.7cm of 1] (3) {};

\node[hidnode,rounded corners] [below =  0.5cm of 1] (4) {};


        \path[color=black,-]
(1) edge  node{} (4)
(2) edge node {} (4)
(3) edge node {} (4);

 \end{tikzpicture}
\end{minipage}
    &
\begin{minipage}[t]{0.45\hsize}
\centering
  \begin{tikzpicture}[baseline=(1), ->,>=triangle 45, shorten >=0pt,
        auto,thick,
        obsnode/.style={circle,inner sep=0.05cm,fill=black,draw,font=\bf, text= white}, 
 hidnode/.style={circle,inner sep=2pt,fill=gray!20,draw,font=\sffamily},
    label/.style={font=\bf}]
        \node[obsnode,rounded corners] (1) {};

\node[obsnode,rounded corners] [below right=  0.35cm and 0.7cm of 1] (2) {};

\node[obsnode,rounded corners] [below right=  0.35cm and 0.7cm of 2] (3) {};


        \path[color=black,-]
(1) edge node {} (2)
(3) edge node {} (2)
;

 \end{tikzpicture}
\end{minipage}\\
    (a) & (b)
  \end{tabular}
  \caption{The possible restrictions of a latent tree to three distinct
  observed variables.  Observed variables correspond
  to solid black dots, latent variables to grey circles.  }
 \label{fig:2tripods}
\end{figure}

\lemref{pseudo} modifies Corollary $1$ in \citet{Shiers} by requiring
the extra equality constraints in \eqref{T3PC} concerning three
distinct variable indices.  For any subset $S \subset {\bf X}$, let
$T|S$ be the restriction of $T$ to $S$, that is, the minimal subtree
of $T$ induced by the elements in $S$ with all the nodes of degree two
not in $S$ suppressed \citep[p.110]{phylo}; refer to \secref{Pfpsuedo} in our supplementary material for the related graphical notions.  \citet{Shiers} only consider
phylogenetic trees in which the observed variables ${\bf X}$ always
correspond to the set of nodes in $T$ with degree one.  In this case
the constraint in \eqref{T3PC} is vacuous.  Indeed, if $X_p, X_q, X_r$
are any three observed nodes in $T$, then $T|\{X_p, X_q, X_r\}$ must
have the configuration on the left panel of \figref{2tripods}, and it
can be seen that
$ph_T(\pi_p , \pi_q) \cup ph_T(\pi_q , \pi_r) \not = ph_T(\pi_p ,
\pi_r) $ for any permutation $(\pi_p, \pi_q,\pi_r)$ of $(p, q, r)$.
However, for a general latent tree $T$ whose observed nodes are
not confined to be the leaves, condition \eqref{T3PC} is necessary for
a pseudo-metric $\delta$ to be $T$-induced: $T|\{X_p, X_q, X_r\}$ may
take the configuration on the right panel of \figref{2tripods}, where
for some permutation $(\pi_p, \pi_q,\pi_r)$ of $(p, q, r)$,
$ph_T(\pi_p , \pi_r) = ph_T(\pi_p , \pi_q) \cup ph_T(\pi_q , \pi_r) $,
and it must hold that
\[
\delta_{\pi_p\pi_r} = \delta_{\pi_p\pi_q} + \delta_{\pi_q\pi_r}
\]
  if $\delta$ is $T$-induced. 

\begin{figure}[t] 
\centering
\begin{minipage}[t]{0.23\hsize}
\centering
 \begin{tikzpicture}[baseline=(1),->,>=triangle 45, shorten >=0pt,
        auto,thick,
        obs node/.style={circle,inner sep=1pt,fill=black,draw,font=\sffamily}, 
  hidden node/.style={circle,inner sep=1.5pt,fill=gray!20,draw,font=\sffamily}]
 \node[obs node,rounded corners] (1) {};
 \node[obs node,rounded corners] [right = 1.2cm of 1] (2) {};
\node[obs node,rounded corners] [below = 1cm of 1] (3) {};
\node[obs node,rounded corners] [below = 1cm of 2] (4) {};
\node[hidden node,rounded corners] [below right= 0.5cm and 0.2cm of 1] (5) {};
\node[hidden node,rounded corners] [below left= 0.5cm and 0.2cm of 2] (6) {};
      \path[color=black!20!black,-,every
      node/.style={font=\sffamily\small}]
  (1) edge node {} (5)
  (3) edge node {} (5)
  (2) edge node {} (6)
  (4) edge node {} (6)
(5) edge node {} (6);

 \end{tikzpicture}
\end{minipage}
\begin{minipage}[t]{0.15\hsize}
\centering
 \begin{tikzpicture}[baseline=(1),->,>=triangle 45, shorten >=0pt,
        auto,thick,
        obs node/.style={circle,inner sep=1pt,fill=black,draw,font=\sffamily}, 
  hidden node/.style={circle,inner sep=1.5pt,fill=gray!20,draw,font=\sffamily}]
 \node[obs node,rounded corners] (1) {};
\node[obs node,rounded corners] [below = 1cm of 1] (2) {};
\node[obs node,rounded corners] [right = 1.3cm of 2] (3) {};
\node[hidden node,rounded corners] [below right= 0.5cm and 0.4cm of 1] (4) {};
\node[obs node,rounded corners] [right= 0.4cm of 4] (5) {};
      \path[color=black!20!black,-,every
      node/.style={font=\sffamily\small}]
(1) edge node {} (4)
 (2) edge node {} (4)
 (3) edge node {} (5)
(4) edge node {} (5);

;
 \end{tikzpicture}
\end{minipage}
 \begin{minipage}[t]{0.2\hsize}
\centering
 \begin{tikzpicture}[baseline=(1),->,>=triangle 45, shorten >=0pt,
        auto,thick,
        obs node/.style={circle,inner sep=1pt,fill=black,draw,font=\sffamily}, 
  hidden node/.style={circle,inner sep=1.5pt,fill=gray!20,draw,font=\sffamily}]
 \node[obs node,rounded corners] (1) {};
 \node[obs node,rounded corners] [below right= 0.25cm and 0.25 cm of 1] (2) {};
\node[obs node,rounded corners] [below right= 0.25cm and 0.25 cm of 2] (3) {};
\node[obs node,rounded corners] [below right= 0.25cm and 0.25 cm of 3] (4) {};
     \path[color=black!20!black,-,every
     node/.style={font=\sffamily\small}]
 (1) edge node {} (4);
\end{tikzpicture}
\end{minipage}
  \begin{minipage}[t]{0.15\hsize}
\centering
 \begin{tikzpicture}[baseline=(1),->,>=triangle 45, shorten >=0pt,
        auto,thick,
        obs node/.style={circle,inner sep=1pt,fill=black,draw,font=\sffamily}, 
  hidden node/.style={circle,inner sep=1.5pt,fill=gray!20,draw,font=\sffamily}]
 \node[obs node,rounded corners] (1) {};
 \node[obs node,rounded corners] [right = 1cm of 1] (2) {};
\node[obs node,rounded corners] [below = 1cm of 1] (3) {};
\node[obs node,rounded corners] [below = 1cm of 2] (4) {};
\node[hidden node,rounded corners] [below right= 0.5cm and 0.5cm of 1] (5) {};
     \path[color=black!20!black,-,every
     node/.style={font=\sffamily\small}]
 (1) edge node {} (5)
 (3) edge node {} (5)
 (2) edge node {} (5)
(4) edge node {} (5);
\end{tikzpicture}
\end{minipage}
 \begin{minipage}[t]{0.19\hsize}
\centering
 \begin{tikzpicture}[baseline=(1),->,>=triangle 45, shorten >=0pt,
        auto,thick,
        obs node/.style={circle,inner sep=1pt,fill=black,draw,font=\sffamily}, 
  hidden node/.style={circle,inner sep=1.5pt,fill=gray!20,draw,font=\sffamily}]
 \node[obs node,rounded corners] (1) {};
 \node[obs node,rounded corners] [below = 0.35cmof 1] (2) {};
\node[obs node,rounded corners] [below right= 0.35cm and 0.35 cm of 2] (3) {};
\node[obs node,rounded corners] [below left= 0.35cm and 0.35 cm of 2] (4) {};
     \path[color=black!20!black,-,every
     node/.style={font=\sffamily\small}]
 (1) edge node {} (2)
 (2) edge node {} (3)
 (2) edge node {} (4);
\end{tikzpicture}
\end{minipage}
\caption{The possible restrictions of a latent tree to four distinct
  observed variables. From left to right, (a)-(e). Observed variables
  correspond to solid black dots, latent variables to grey circles.
}
 \label{fig:config}
\end{figure}
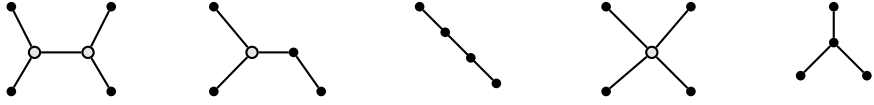

  While condition \eqref{T4PC} appears in the result of
  \citet{Shiers}, it may lead to different patterns of constraints for a general latent tree. For four
  distinct indices $1 \leq p, q, r, s \leq m$, there are three
  possible partitions into two subsets of equal sizes, namely,
  $\{p, q\} |\{ r, s\}$, $\{p, r\} | \{q, s\}$ and
  $\{p, s\} | \{q, r\}$.  These three partitions correspond to the path pairs
  \begin{equation} \label{pathPairs}
  (ph_T(p, q), ph_T(r, s)), (ph_T(p, r), ph_T(q, s)) \text{ and }(ph_T(p, s), ph_T(q, r))
\end{equation}
respectively. Now refer to \figref{config} which shows all possible
configurations of the restriction of $T$ to the four observed
variables $X_p, X_q, X_r, X_s$.  In \figref{config}(a)-(c), up to
permutations of the indices $\{p, q, r, s\}$, only one of three pairs
in \eqref{pathPairs} can give an empty set when the intersection of
its two component paths is taken.  In light of \eqref{T4PC}, this
implies that, for some permutation $\pi$ of the indices $p, q, r, s$,
\begin{equation} \label{permutationineq}
\delta_{\pi_p \pi_q} + \delta_{\pi_r \pi_s } \leq \delta_{\pi_p \pi_r} + \delta_{\pi_q \pi_s }  =\delta_{\pi_p \pi_s} + \delta_{\pi_q\pi_r }. 
\end{equation}
%
On the contrary, in \figref{config}(d) and (e), it must be the case
that each of the three path pairs in \eqref{pathPairs} gives an empty
set when an intersection is taken between its two component paths,
giving the equalities
$
\delta_{pq} + \delta_{rs } =  \delta_{p r} + \delta_{q s }  =\delta_{p s} + \delta_{q r }
$
in consideration of \eqref{T4PC}. 

\lemref{pseudo} readily implies a characterization of the latent tree model
$\mathcal{M}_{\bf X}(T)$ via polynomial constraints in the entries of
the covariance matrix $\Sigma=(\sigma_{pq})$ as spelt out in the ensuing corollary. Its proof employs similar arguments in \citet{Shiers} and is deferred to our supplementary material. In what follows, we let
$\mathcal{Q} \subset {m \brace 4}$ be the set of all quadruples
$\{p, q, r, s\} \in {m \brace 4}$ such that only one of the three path
pairs in \eqref{pathPairs} gives an empty set when the union of its
two component paths is taken.  In other words, $\mathcal{Q}$ contains
all $S \in {m \brace 4}$ such that $T|S$ is one of the configurations
in \figref{config}(a)-(c).  Given $\{p, q, r, s\} \in \mathcal{Q}$, we
write $\{p, q\}|\{r, s\} \in \mathcal{Q}$ to indicate that
$\{p, q, r, s\}$ belongs to $\mathcal{Q}$ in a way that it is the path
pairs $ph_T(p , q)$ and $ph_T(r , s)$ that have empty intersection.
Similarly, we will let $\mathcal{L}$ be the set of all triples
$S = \{p , q, r\}\in {m \brace 3}$ such that $T|S$ has the
configuration in \figref{2tripods}$(b)$.   We will use the notation
$p - q- r \in \mathcal{L}$ to indicate that $q$ is the ``middle point"
such that $ph_T(p , q) \cap ph_T(q , r) = \emptyset$.

\begin{corollary} \label{cor:charGT} Suppose
  $\Sigma = (\sigma_{pq})_{1 \leq p , q \leq m}$ is the covariance
  matrix of $\bf X$ and has no zero entries. The following are
  together necessary and sufficient for the distribution of $\bf X$ to
  belong to $\mathcal{M}_{\bf X}(T)$:
\begin{enumerate}
\item Inequality constraints: 
\begin{enumerate}
\item For any $\{p, q, r\} \in {m \brace 3}$,  $\;\sigma_{pq} \sigma_{pr} \sigma_{qr} \geq 0. $
\item For any $\{p, q, r\} \in {m \brace 3}\backslash  \mathcal{L}$, 
\[
\sigma_{pq}^2 \sigma_{qr}^2 - \sigma_{qq}^2\sigma_{pr}^2, \ \ \sigma_{pr}^2 \sigma_{qr}^2 - \sigma_{rr}^2\sigma_{pq}^2, \ \ \sigma_{pq}^2 \sigma_{pr}^2 - \sigma_{pp}^2 \sigma_{qr}^2  \leq 0.\]
%
\item For any $\{p, q\}|\{r, s\} \in \mathcal{Q}$,  $\sigma_{pr}^2\sigma_{qs}^2 - \sigma_{pq}^2\sigma_{rs}^2  \leq 0.$
\end{enumerate}
\item Equality constraints: 
\begin{enumerate}
\item For any $p - q - r \in \mathcal{L}$, $\sigma_{pq}\sigma_{qr} - \sigma_{qq} \sigma_{pr} = 0$. 
\item For any $\{p, q\}|\{r, s\}  \in \mathcal{Q}$, $\sigma_{pr}\sigma_{qs} - \sigma_{ps} \sigma_{qr} = 0.$
\item For any $\{p, q, r, s\} \not \in \mathcal{Q}$, 
$
\sigma_{ps} \sigma_{qr} - \sigma_{pr}\sigma_{qs}  = \sigma_{pq}\sigma_{rs} - \sigma_{pr}\sigma_{qs}= 0.
$
\end{enumerate}
\end{enumerate}
\end{corollary}

\section{Testing a star tree model} \label{sec:test}

In this section we illustrate how one can test a postulated Gaussian
latent tree model using \corref{charGT}.  In order to focus the discussion we treat the
simple but important special case of a star tree, which corresponds to a single factor model. 
A single factor model with $m$ observed variables
${\bf X} = \{X_1\dots, X_m\}$ can  be described by the linear
system of equations
\begin{equation} \label{1factor}
X_p =  \mu_p + \beta_p H + \epsilon_p, \qquad 1 \leq p \leq m,
\end{equation}
where $\mu_p$ is the mean of $X_p$, $H \sim N(0, 1)$ is a latent
variable, $\beta_p$ is the loading coefficient for variable $X_p$, and
$\epsilon_p \sim N(0, \sigma_{p, \epsilon}^2)$  is the idiosyncratic error for variable
$X_p$.  All of $H$, $\epsilon_1, \dots , \epsilon_m$ are independent. The model postulates that $X_1,\dots,X_m$ are
conditionally independent given $H$.  It thus corresponds to the
graphical model associated with a star tree $T_\star = (V, E)$ with
$V ={\bf X}\cup\{ H\}$, $
E = \{ (H, X_p)\}_{ 1 \leq p \leq m}.$ 

Let ${\bf X}_1,\dots,{\bf X}_n$ be i.i.d.~draws from the distribution
of $\bf X$, which is assumed to be Gaussian.  Our goal is to test
whether the distribution of $\bf X$ belongs to the single factor model
$\mathcal{M}_{\bf X}(T_\star)$.  Without loss of generality, we may
assume that $\mu_p = 0$ for all $p \in [m]$ \citep[Theorem
3.3.2]{Anderson}.  We proceed by testing whether all the constraints
in \corref{charGT} are simultaneously satisfied with respect to the
latent tree $T_\star$.  For simplicity, we will focus on testing the
equality constraints in \corref{charGT}$(ii)$, and briefly discuss how
one can incorporate the inequality constraints in \corref{charGT}$(i)$
in \secref{conclude}.  For $T_\star$, both sets $\mathcal{L}$ and
$\mathcal{Q}$ are empty, so that \corref{charGT}$(ii)(a)$ and $(b)$
are automatically satisfied.  Hence, we are only left with
\corref{charGT}$(ii)(c)$:
For any $\{p, q, r, s\} \in {m \brace 4}$, 
\begin{equation} \label{tetrad}
 \sigma_{ps} \sigma_{qr} - \sigma_{pr}\sigma_{qs}  = \sigma_{pq}\sigma_{rs} - \sigma_{pr}\sigma_{qs}= 0.
\end{equation}
The two polynomials above, equal to $\text{det} (\Sigma_{pq, sr})$ and
$\text{det} (\Sigma_{ps, qr})$ respectively,  are known as
\emph{tetrads} in the literature of factor analysis. It is well-known
that they define equality constraints for a single factor model
\citep{deLeeuw, bollen1993confirmatory,algfac}.


\subsection{Estimating tetrads} \label{sec:estimate-tetrads}

The idea now is to estimate each one of the $2 \cdot {m \choose 4}$
tetrads in \eqref{tetrad}, and aggregate the estimates in a test
statistic.  From the sample covariance matrix 
  $S = (s_{pq}) = n^{-1} \sum_{i = 1}^n {\bf X}_i {\bf X}_i^T$, a straightforward sample tetrad estimate, say $ s_{ps} s_{qr} - s_{pr}s_{qs}$, can be computed.  If one define the vectors ${\bf t}  = (s_{ps}, s_{qr} , s_{pr},  s_{qs})'$ and ${\bf t}_0 = (\sigma_{ps}, \sigma_{qr} , \sigma_{pr},  \sigma_{qs})'$, as well as the function $g({\bf t}) = s_{ps} s_{qr} - s_{pr} s_{qs}$, by the delta method it is  expected that $\sqrt{n} (g({\bf t}) - g({\bf t}_0)) \rightarrow N(0, 	\nabla g({\bf t}_0)'V \nabla g({\bf t}_0))$, where $V$ is the limiting covariance matrix of $\sqrt{n}( {\bf t} - {\bf t}_0)$ and $\nabla g({\bf t}_0)$ is the gradient of $g(\cdot)$ evaluated at ${\bf t}_0$.  
 However, the distribution of
this sample tetrad becomes asymptotically degenerate at singularities,
that is, when the gradient $\nabla g({\bf t}_0)$ vanishes, which happens if the
underlying true covariances are zero 
\citep{drton:xiao:2016}.  Consequently, a standardized sample tetrad
cannot be well approximated by a normal distribution if the underlying
correlations are weak.   More generally, even for stronger
correlations, we found it difficult to reliably estimate the variance
of all sample tetrads in larger-scale models. 

We propose alternative estimators for which
sampling variability can be estimated more easily.   Due to the independence of
samples, the tetrad
$\text{det} (\Sigma_{pq, sr}) = \sigma_{ps} \sigma_{qr} -
\sigma_{pr}\sigma_{qs}$ can be estimated {unbiasedly} with the
differences
\begin{equation}
  \label{eq:simple-estimates}
Y_{i, (pq)(sr)}:= 
X_{p, i} X_{s, i} X_{q, i+1} X_{r, i+1} - X_{p, i} X_{r, i} X_{q, i+1} X_{s, i+1}, \ \  i = 1, \dots, n-1, 
\end{equation}
where the subscripts in $Y_{i, (pq)(sr)}$ is indicative of the row and
column indices for the submatrix $\Sigma_{pq, sr}$.  These differences
can then be averaged for an estimate of the tetrad.  Similarly, one
can form $Y_{i, (ps)(qr)}$ to estimate $\text{det} (\Sigma_{ps, qr})$
in \eqref{tetrad}. If we arrange all the tetrads from
$\{\text{det} (\Sigma_{pq, sr}) , \text{det} (\Sigma_{ps, qr})\}_{\{p,
  q, r, s\} \in {m \brace 4}}$ into a $2  {m \choose 4}$-vector
$\Theta$, and correspondingly arrange the estimates
$\{Y_{i, (pq)(sr)}, Y_{i, (ps)(qr)}\}_{\{p, q, r, s\} \in {m \brace
    4}}$ into a $2 {m \choose 4}$-vector ${\bf Y}_i$ for each
$i$, then the central limit theorem for \emph{1-dependent} sums
ensures that for sufficiently large sample size $n$ we have the
distributional approximation
\begin{equation}\label{CLT}
\sqrt{n-1} (\bar{\bf Y}  - \Theta) \approx_d N(0,  \Upsilon),
\end{equation}
where $\bar{\bf Y}= (n-1)^{-1} \sum_{i = 1}^{n-1} {\bf Y}_i$ and 
$
\Upsilon = \text{Cov}[{\bf Y}_1, {\bf Y}_1] + 2 \text{Cov}[{\bf Y}_1, {\bf Y}_2]$. The latter limiting covariance matrix will not degenerate to a singular matrix even if the underlying covariance matrix for ${\bf X}$ has zeros at which some of the tetrads are singular (i.e. have zero gradient).

\subsection{Bootstrap test} \label{sec:estimate-tetrads}

The fact from~\eqref{CLT} could serve as the starting point for a test
of model $\mathcal{M}(T_\star)$.  However, the normal approximation
quickly becomes of concern when moving beyond a small number of
variables $m$.  Indeed, the dimension of $\Theta$, $2 {m \choose 4}$,
may well be close to the sample size $n$, or even larger.  For
instance, if $n = 250$, for a model with merely $8$ observed variables
the dimension of $\Theta$ is already $ 2 {8 \choose 4} = 140$, more
than half the sample size. A recent work of \citet{WuZhang}, which
follows up on the groundbreaking paper of \citet{CCKHD} on Gaussian
approximation for maxima of high dimensional independent sums,
suggests that while the approximation in \eqref{CLT} may be dubious,
by taking a supremum norm on both sides, the Gaussian approximation
\begin{equation} \label{GaussainApprox}
\sqrt{n-1} \|(\bar{\bf Y}  - \Theta)\|_\infty \approx_d \|{\bf Z}\|_\infty, 
\end{equation}
where ${\bf Z} =_d N(0, \Upsilon)$, can be valid even the dimension of
$\Theta$ is large compared to $n$. In fact, the original work of
\citet{CCKHD} suggested that asymptotically, the dimension can be
sub-exponential in the sample size for the Gaussian approximation to
hold.  In what follows, we will discuss implementation of
and experiments with a vanishing tetrad test based on
\eqref{GaussainApprox}.  While it is possible to adapt the supporting
theory for the present application, the technical details are involved
and beyond the scope of this
conference paper.

Since $\bar{\bf Y}$ from~(\ref{CLT}) and~(\ref{GaussainApprox}) is an
estimator of the vector of tetrads $\Theta$, it is natural to use
$\|\bar{{\bf Y}}\|_\infty$ as the test statistic and reject the model
$\mathcal{M}(T_\star)$ for large values of $\|\bar{\bf
  Y}\|_\infty$. The Gaussian approximation \eqref{GaussainApprox}
suggests that when $\mathcal{M}(T_\star)$ is true, i.e. $\Theta = 0$,
$\sqrt{n-1}\|{\bf Y}\|_\infty$ is distributed as $\|{\bf
  Z}\|_\infty$. Nevertheless, to calibrate critical values based on
the distribution of $\|{\bf Z}\|_\infty$, one must estimate the
unknown covariance matrix $\Upsilon$.  \citet{WuZhang} suggested  the batched mean estimator
\begin{equation} \label{hatUpsilon}
\hat{\Upsilon} = \frac{1}{B \omega}\sum_{b = 1}^{\omega} \left(\sum_{i \in L_b} ({\bf Y}_i - \bar{{\bf Y}})\right) \left(\sum_{i \in L_b} ({\bf Y}_i - \bar{{\bf Y}})\right)^T,
\end{equation}
where for a batch size $B$ and $\omega := \lfloor (n-1)/B\rfloor$ one
considers the non-overlapping sets of samples
$L_b = \{1 + (b-1)B, \dots, bB\}$, $b = 1, \dots, \omega$.  The
``batching" aims to capture the dependence among the ${\bf Y}_i$'s,
and has been widely studied in the time series literature
\citep{buhlmann2002bootstraps, Lahiri}.  If model
$\mathcal{M}(T_\star)$ is true, then 
\eqref{GaussainApprox} yields that 
\begin{equation*}
\mathcal{T} := \sqrt{n-1} \|\diag(\hat{\Upsilon})^{-1/2}\bar{\bf Y} \|_\infty  \approx_d 
\|\diag(\hat{\Upsilon})^{-1/2}\tilde{\bf Z}\|_\infty,
\end{equation*}
where the right-hand side is to interpreted conditionally on
$\hat{\Upsilon}$, with $\tilde{\bf Z} \sim N(0, \hat{\Upsilon})$ and
$\diag(\hat{\Upsilon})$ comprising only the diagonal of
$\hat{\Upsilon}$.  More precisely, for a fixed test level
$\alpha \in (0,1)$, if we define $q_{1- \alpha}$ to be the
\emph{conditional} $(1- \alpha)$-quantile of the distribution of
$\|\diag(\hat{\Upsilon})^{-1/2}\tilde{\bf Z}\|_\infty$ given
$\hat{\Upsilon}$, then
\begin{equation} \label{cor5.4}
P(\mathcal{T} > q_{1-\alpha}) \approx \alpha,
\end{equation}
according to \citet[Corollary 5.4]{WuZhang}.  We will use
$\mathcal{T}$ as our test statistic for the model
$\mathcal{M}(T_\star)$, and calibrate the critical value based on
\eqref{cor5.4}  by simulating the conditional quantile $q_{1 -
  \alpha}$ from $\|\diag(\hat{\Upsilon})^{-1/2}\tilde{\bf Z}\|_\infty$
for fixed $\hat{\Upsilon}$.


\subsection{Implementation }

While our above presentation invoked the estimate $\hat\Upsilon$,
which is a matrix with $O(m^8)$ entries, we may in fact bypass the
problem of computing such a large covariance matrix for the tetrad
estimates.  To simulate the conditional quantile $q_{1- \alpha}$ in
\eqref{cor5.4}, let $e_1, \dots, e_{\omega}$ be
i.i.d.~standard normal random variables, and consider the expression
\begin{equation} \label{expression}
\left\|\frac{\diag(\hat{\Upsilon})^{-1/2}}{\sqrt{B \omega}} \sum_{b = 1}^{\omega} e_b\left(\sum_{i \in L_b} ({\bf Y}_i - \bar{{\bf Y}})\right)\right\|_\infty, 
\end{equation}
which has exactly the same distribution as
$\|\diag(\hat{\Upsilon})^{-1/2}\tilde{\bf Z}\|_\infty$ conditioning on
the data ${\bf X}_1, \dots, {\bf X}_n$.  We emphasize the $O(m^4)$
diagonal entries of $\hat{\Upsilon}$ are easily computed as variances
in~(\ref{hatUpsilon}).  In conclusion, we perform the following
multiplier bootstrap procedure:
\begin{inparaenum}[(i)]
\item Generate many,
  say $E = 1000$, sets of $\{e_1, \dots, e_{\omega}\}$,
\item evaluate
  \eqref{expression} for each of these $E$ sets, and
\item take
  $q_{1- \alpha}$ to be the $1 - \alpha$ quantile from the resulting
  $E$ numbers.
\end{inparaenum}
Despite the bootstrap being a computationally intensive process, it is
not hard to see that the evaluation of \eqref{expression} for all $E$
sets of multipliers will involve $O(m^4 n E)$ operations, which
even for moderate $m$ is far less than the $O(m^8)$ operations needed
to obtain an entire covariance matrix for all tetrads.

\begin{remark}
  It is instructive to make a comparison with the testing methodology
  in \citet{Shiers}, where the focus was on lower-dimensional
  applications.  Suppose $\tau: \Sigma \mapsto \Theta$ is the function
  that maps the covariance matrix $\Sigma$ into the vector $\Theta$ of
  tetrads in \eqref{tetrad}.  To test the vanishing of the tetrads,
  \citet{Shiers} form plug-in estimates $\hat{\Theta} = \tau(S)$ for
  $\Theta$ with the sample covariance matrix
  $S = n^{-1} \sum_{i = 1}^n {\bf X}_i {\bf X}_i^T$.  Letting $\text{Var}[\tau(S)]$ be the covariance matrix for the
  $2{m \choose 4}$-vector $\tau(S)$, they form a Hotelling's $T^2$
  type statistic as
\begin{equation} \label{Hotelling}
n \tau(S)^T (\widehat{\text{Var}}[\tau(S)])^{-1} \tau(S),
\end{equation}
where $\widehat{\text{Var}}[\tau(S)]$ is a consistent estimate for
$\text{Var}[\tau(S)]$; see also \cite{MorMinors}.  For a test of model
$\mathcal{M}(T_\star)$, this statistic is now compared to a chi-square
distribution with $2 {m \choose 4}$ degrees of freedom.  While this
calibration is justified for sufficiently large sample size $n$ by a
joint normal approximation analogous to~(\ref{GaussainApprox}), it can
be problematic for large  $m$.  Even more pressing
can be the computational disadvantage that one explicitly uses the
entire matrix $\widehat{\text{Var}}[\tau(S)]$ with its $O(m^8)$ entries.
\end{remark}

\section{Numerical experiments}
\label{sec:num-exp}

We now report on some experiments with the bootstrap test based on the
sup-norm of the estimated tetrads $\mathcal{T}$ proposed in \secref{test}. In the implementation we always use $E = 1000$ sets of normal multipliers
to simulate the quantile $q_{1 - \alpha}$ and work with batch size
$B = 3$ in \eqref{expression}. We also
benchmark our methodology against the likelihood ratio test for factor
models implemented by the function \texttt{factanal} in the
\texttt{base} library of \texttt{R}, which implements a likelihood
ratio (LR) test with Bartlett correction for more accurate asymptotic
approximation. The critical value of the LR test is
calibrated with the chi-square distribution with
${m-1 \choose 2} - 1$ degrees of freedom
\citep[p.99]{MR2723140}. 

\subsection{Low dimensional setup}
We first consider two experimental setups, each with data generated from the one-factor model in \eqref{1factor} for  both $(m, n) = (20, 250)$ and $(m, n )= (20, 500)$.  The model parameters are as follows: \begin{inparaenum} [(i)]\item Setup 1: all loadings $\beta_p$ and error variances $\sigma^2_{p, \epsilon}$ are taken to be $1$. \item Setup 2: $\beta_1$ and $\beta_2$ are taken to be $10$, while the other loadings are independently generated based on a normal distribution with mean $0$ and variance $0.2$. The error variances $\sigma^2_{p, \epsilon}$ all equal $1/3$. \end{inparaenum}

\begin{figure}[t]
\centering
\begin{minipage}[t]{.4\linewidth}
 \includegraphics[width=1\linewidth]{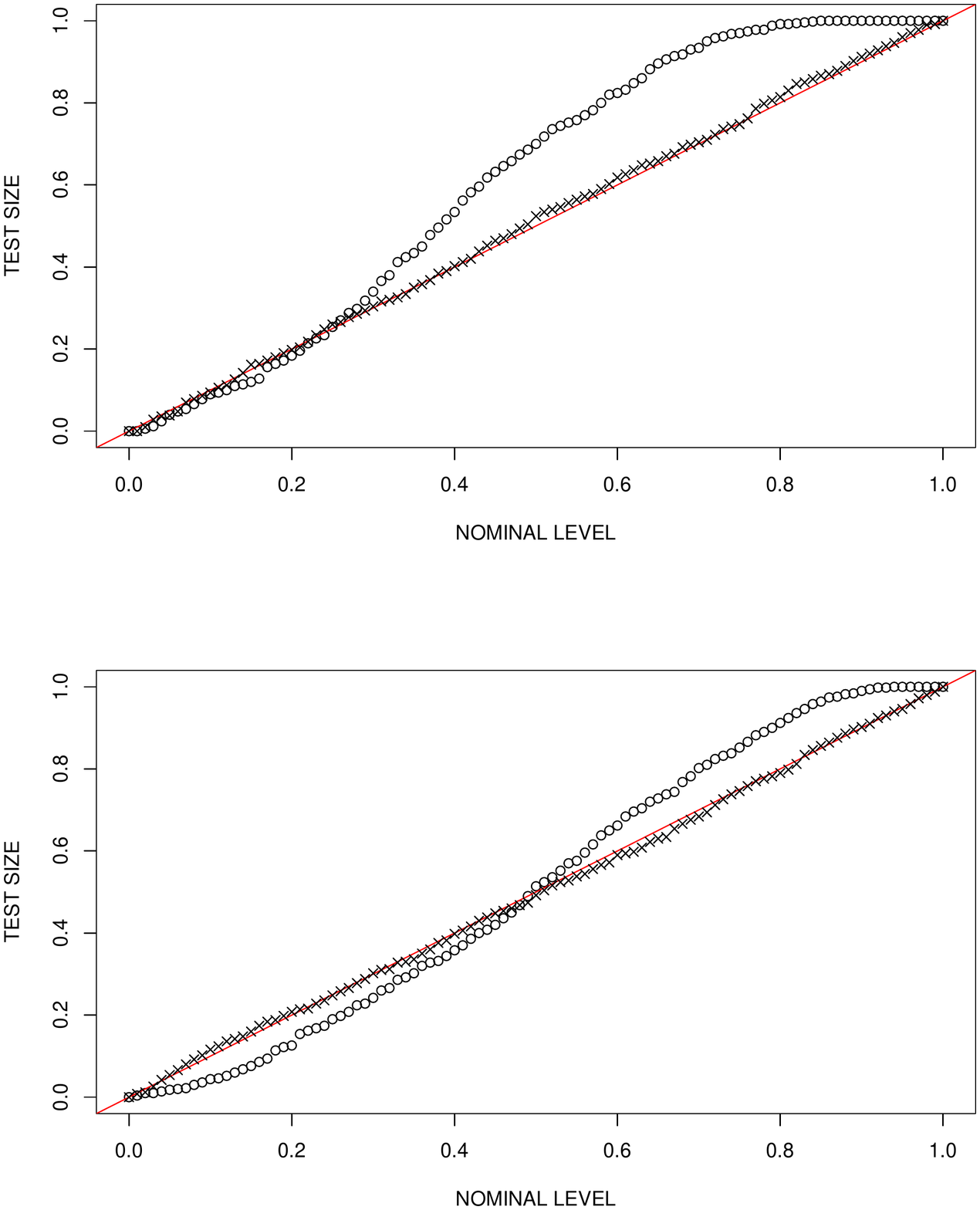}
 \end{minipage} 
 \begin{minipage}[t]{.4\linewidth}
 \includegraphics[width=1\linewidth]{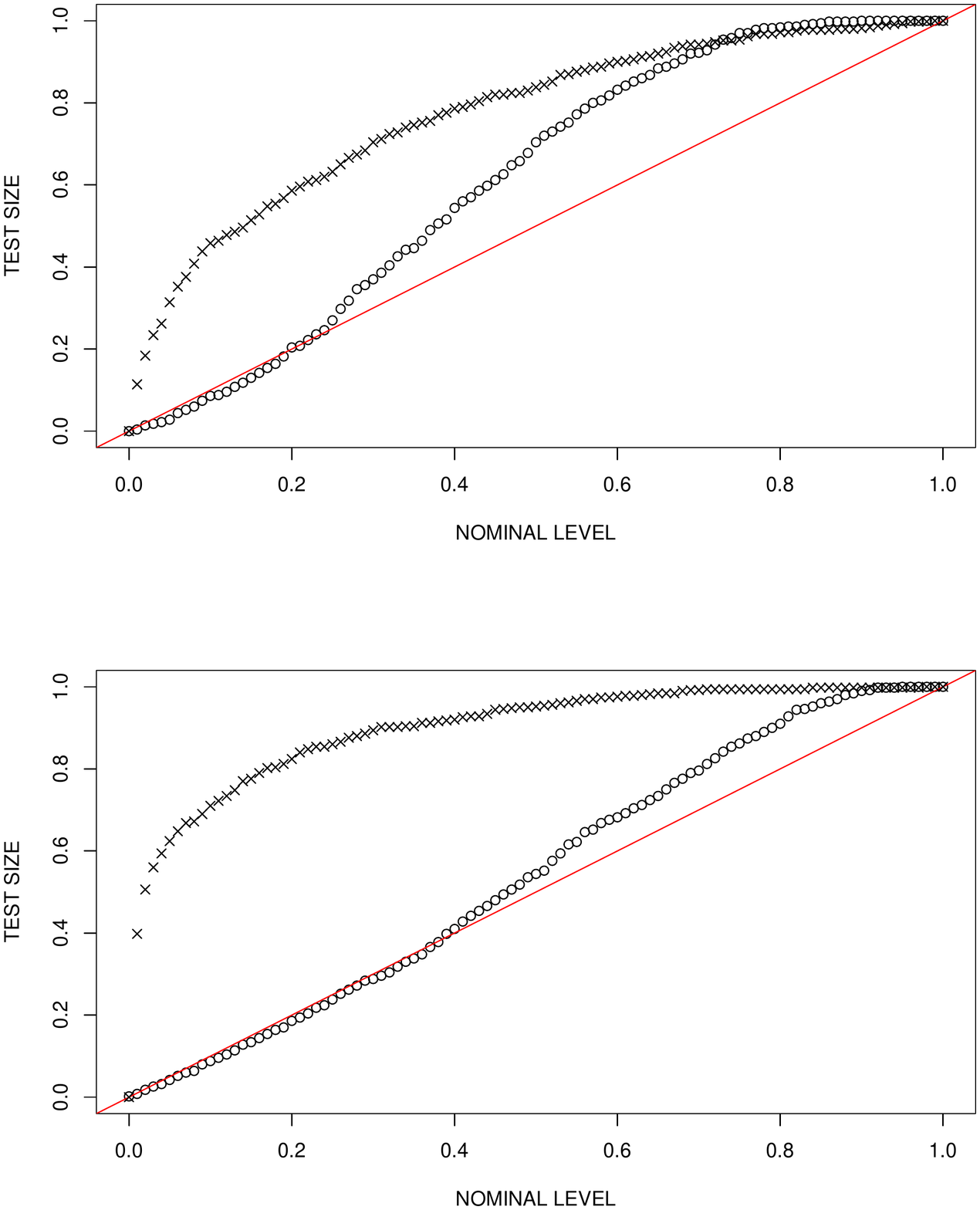}
 \end{minipage}
 \caption{Empirical test sizes vs nominal test levels based on $500$ experiments. Data are generated based on $\mathcal{M}_{\bf X}(T_\star)$ with parameters as prescribed in the text. Upper panels: $(m,n) = (20, 250)$. Lower panels: $(m,n) = (20, 500)$. Left panels: Setup 1. Right panels: Setup 2. Open circles: Test based on the statistic $\mathcal{T}$.  Crosses: LR test implemented by \texttt{factanal}.}
 \label{fig:LDfig} 
\end{figure}

For different nominal test levels $\alpha$ in the range $(0, 1)$ that are $0.01$
apart, we compare the empirical sizes of our test based on the
statistic $\mathcal{T}$ and the likelihood ratio (LR) test implemented
by the function \texttt{factanal}, using $500$
repetitions of experiments.  
  The results are shown in \figref{LDfig}. The left two panels correspond to Setup 1 and the right two panels correspond to Setup 2, while the upper panels correspond to $(m, n) = (20, 250)$ and lower correspond to $(m, n )= (20, 500)$. 
While we show the entire range $(0,1)$ for the x-axis, practical interest is typically in the initial part where the nominal error rate is in say $(0,0.1)$. 

\begin{figure}[t]
\centering
 \includegraphics[height=0.4\linewidth, width = 0.5 \linewidth]{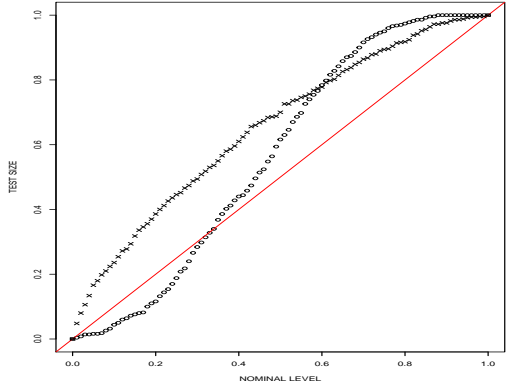}
  \caption{Empirical test size vs nominal test levels based on $500$ experiments for data generated from $\mathcal{M}_{\bf X}(T_\star)$ under Setup 1 and $(m, n) = (100, 250)$. Open circles: Test based on $\mathcal{T}$.  Crosses: LR test implemented by \texttt{factanal}. }
  \label{fig:whtever3}   
\end{figure}

In Setup $1$, for both sample sizes, the empirical test sizes of the LR test align almost perfectly with the $45^\circ$ line as one would expect from classical theory. The sizes of our test based on $\mathcal{T}$ also align better with $45^\circ$ line as sample sizes grow. Note that for nominal test levels that are of practical interest, $\mathcal{T}$ also gives conservative test sizes for both sample sizes.



In Setup $2$, where parameters are close to being ``singular", one can see the true advantage of using $\mathcal{T}$ over the LR test. The empirical test sizes of the LR test with
\texttt{factanal} do not align well with the $45^\circ$ line as
one normally expect from classical theory, whereas the test sizes of our statistic
$\mathcal{T} $ lean closer to the $45^\circ$ line as $n$
increases. Particularly  the performance of the LR test is  problematic since, by rejecting the true model \eqref{1factor} all too often,  it fails to give even an approximate control on type $1$ error. 
Note that the values of $\beta$ and $\sigma_{p, \epsilon}$
are such that, for the most part, the observed variables $\bf X$ are rather weakly
dependent on each other.  If the observations were in fact
independent then the likelihood ratio test statistic does not exhibit
a chi-square limiting distribution \citep[Theorem 6.1]{drton:lrt}.  This highlights the fact that, in addition to avoiding
any non-convex optimization of the likelihood function of the factor
model, our approach based on the simple estimates from~(\ref{eq:simple-estimates}) is not
subject to non-standard limiting behaviors that plague the LR test when the parameter values lean close to singularities of the parameter space \citep{drton:lrt}.

\subsection{Higher dimensional setup}

Our last experiment aims to compare the test sizes of the two tests when the number of observed variables $m$ is relatively large compared to $n$. Data are exactly as in Setup $1$, except that $(m, n) = (100, 250)$.  For such a model with large $m$,  the number of tetrads involved in our testing methodology is so large that even after taking the supremum norm  one shouldn't expect \eqref{GaussainApprox}  to hold; for example, when $m = 50$, the dimension of $\Theta$ is $2 \cdot {50 \choose 4} = 460600$, and one should be skeptical about the validity of \eqref{GaussainApprox}  when we only have the sample size $n = 250$. To implement our test,  we first randomly select $10000$ of the $2 \cdot {m \choose 4} $ tetrads, and proceed with the bootstrapping procedure in \eqref{expression} with ${\bf Y}_i$ being estimates for this selected subset of tetrads alone. The choice of $10000$ tetrads to be tested is based on the fact that,  in the previous experiments with $(m, n) = (20, 250)$, our test gives reasonable empirical test sizes for a  practical range of nominal levels when the total number of tetrads being tested, $2 \cdot {20 \choose 4}$,  is approximately $10000$. Since the subset of tetrads is randomly selected, our test is still expected to approximately control the test sizes at nominal level. The results are reported in \figref{whtever3} .

As seen, the test based on $\mathcal{T}$ has
the main features seen in  the first experiment.  In particular, it
successfully controls type I error rates for the practical range of
$\alpha\in(0, 0.1)$.  In contrast, with $m$ increased to $100$,
the LR test drastically fails to control type I error rate.  This is
despite the fact that the setup is regular with parameter values that
are far from any model singularity.  The reason for the failure of the
LR test is the fact that the dimension is on the same order as the
sample size of $250$.  The sample size is not large enough for
chi-square asymptotics based on fixed dimension $m$ to ``kick in''.


\section{Discussion} \label{sec:conclude}

In this paper we have established a full  set of polynomial constraints on the covariance matrix of the observed variables, in the form of both equalities and inequalities, that characterizes a general Gaussian latent tree model whose observed nodes are not confined to be the leaves.  Focusing on the special case of a star tree model, we also experimented with a new methodology for testing the equality constraints by forming unbiased estimates of the polynomials involved. In simulation studies, when the number of variables involved is large or the underlying parameters are close to being ``singular", our test compares favorably with the likelihood ratio test in terms of test size. 


Our results have paved the way for developing a full-fledged algebraic test for a Gaussian latent tree model. Although we have not pursued this generality in the present conference paper, we give a brief discussion here. Of course, to do so one would first need to write an efficient graph algorithm to tease out all the polynomials entailed by \corref{charGT} for a given latent tree input. Then the current  testing methodology can be adopted  by forming unbiased estimates of all these polynomials at hand, which also brings to our attention that 
in \secref{test} only the equality constraints in \corref{charGT}$(ii)$ were used to test the single factor model. For illustration, 
take the 3-degree monomial in \corref{charGT} (i)(a) as an example. Like~(\ref{eq:simple-estimates}), one may form a summand
$
Y_{i, (p, q, r)} = X_{p, i} X_{q, i} X_{p, i+1} X_{r, i+1} X_{q, i+2} X_{r, i+2}, 
$
which is unbiased for $\sigma_{pq} \sigma_{pr}\sigma_{qr}$, and then use $(n-2)^{-1}\sum_{i = 1}^{n-2} Y_{i, (p, q, r)} $ as an averaged estimator.  To incorporate the constraints in  \corref{charGT} (i) into our test one can first arrange all those inequalities into ``less than" conditions, i.e., \corref{charGT} (i)(a) becomes $- \sigma_{pq} \sigma_{pr} \sigma_{qr} \leq 0$ and the corresponding estimate becomes $- (n-2)^{-1}\sum_{i = 1}^{n-2} Y_{i, (p, q, r)} $. Following that,  in the definition of the test statistic $\mathcal{T}$, one can take a maximum over all the unbiased estimates for the ``less than" versions of the polynomials in \corref{charGT}$(i)$,  in addition to the absolute values of the estimates for the polynomials in  \corref{charGT}$(ii)$. The resulting test statistic shall also reject the model $\mathcal{M}(T_\star)$ when its value is too large. While critical values can still be calibrated with multiplier bootstrap, additional techniques such as inequality selection can be incorporated to contain the power loss as a result of testing the inequalities; see \citet{CCKmoments} for more details.

Another challenge is the determination of the batch size $B$ in \eqref{hatUpsilon}. In our simulation studies of \secref{num-exp} we took $B = 3$ since we believe that a batch size of $3$ should be enough to capture dependence among the $1$-dependent summands. 
Batch size determination has been widely studied in the time series literature for low dimensional problems \citep{buhlmann2002bootstraps, Hall, Lahiri}. To the best of our knowledge,  in high dimensions this is still a widely open problem. Theoretical research on this is far beyond the scope of our current work.

\section{Supplementary material}

In this supplement we furnish proofs for the main text of ``Algebraic tests of general Gaussian latent tree models". 
\subsection{Proof of \corref{charGT}}
We only sketch the proof here since it is exactly analogous to that of Theorem $3$ in \citet{Shiers}. First, consider the special case where all the entries of $\Sigma$, and hence the the Pearson correlations $\rho_{pq}$, $1 \leq p \not = q \leq m$,  are strictly positive.  In this case condition $(i)(a)$ is redundant. Via the isomorphsim
\[
\delta_{pq} = - \log \rho_{pq}, 
\]
between the parametrizations in \eqref{para} and all $T$-induced
pseudometrics,  the discussion preceding our corollary readily
translates \eqref{T4PC} into $(i)(c), (ii)(b), (ii)(c)$ and
\eqref{T3PC} into $(ii)(a)$, whereas the triangular inequality
property of pseudometrics is translated into $(i)(b)$ for triples
$\{p, q, r\}$ that are not in $\mathcal{L}$.
The general case of $\Sigma$ with nonzero but not necessarily positive
entries is then addressed by incorporating condition $(i)(a)$.

\subsection{Proof of \lemref{pseudo}} \label{sec:Pfpsuedo}

To prove the lemma, we first collect all the required graphical notions borrowed from \citet{phylo}. We attempted to make this proof as self-contained as possible, but the readers are encouraged to read \citet{phylo} for more background on mathematical phylogenetics.

Suppose we are given a tree $T = (V, E)$.  If $\tilde{V}$ is a subset of $V$, $T(\tilde{V})$ denotes the minimal subtree of $T$ that contains all the nodes in $\tilde{V}$. If $e \in E$, $T \backslash e$ is the graph obtained by removing $e$, and $T/e$  is the tree obtained from $T$ by identifying the ends of $e$ and then deleting $e$. In particular, if $v \in V$ is a node of degree two and $e$ is an edge incident with $v$, $T/e$ is said to be obtained from $T$ by \emph{suppressing} $v$. If $v_1, v_2 \in V$, $ph_T(v_1, v_2)$ is the set of edges on the unique path connecting $v_1$ and $v_2$. 

We will also need the notion of an ${\bf X}$-tree. An ${\bf X}$-tree, or semi-labeled tree on a set ${\bf X}$,  is an ordered pair $\mathcal{T} = (T, \phi)$, where $T$ is a tree with node set $V$  and $\phi: {\bf X} \rightarrow V$ is a (labeling) map with the property that, for each $v \in V$ \emph{of degree at most two}, $v \in \phi({\bf X})$. Note that $\phi$ is not necessarily injective. Moreover, if ${\bf X}'$ is a subset of ${\bf X}$, $T|{\bf X}'$ is the tree obtained from $T(\phi({\bf X}'))$ by suppressing all the nodes of degree two that are not in $\phi({\bf X}')$. We then define the \emph{restriction of} $\mathcal{T}$ \emph{to} ${\bf X}'$, denoted $\mathcal{T}|{\bf X}'$, to be the ${\bf X}'$-tree $(T|{\bf X}', \phi |{\bf X}')$.

Finally, we introduce the notion of ${\bf X}$-\emph{split}. For a set ${\bf X}$, an ${\bf X}$-split is a partition of ${\bf X}$ into two non-empty sets. We denote the ${\bf X}$-split whose blocks are $A$ and $B$ by $A|B$ where the order of $A$ and $B$ in the notation doesn't matter. Now suppose $\mathcal{T} = (T, \phi)$ is an ${\bf X}$-tree with an edge set $E$. For each  $e \in E$, $T \backslash e$ must consist of two components $V^e_1$ and $V^e_2$ which induce an ${\bf X}$-split $\phi^{-1}(V^e_1) |\phi^{-1}(V^e_2)$. We then define $\Sigma(\mathcal{T}):= \{ \phi^{-1}(V^e_1) |\phi^{-1}(V^e_2) : e \in E\}$ as the collection of all ${\bf X}$-splits induced by $\mathcal{T}$. 

\emph{Important remark}: In all the definitions above,  ${\bf X}$ is not specified as a subset of the node set $V$ for a given tree. Nonetheless, when we have a tree $T = (V, E)$ with a subset of observed nodes ${\bf X}\subset V$  as in the main text,  we will slightly abuse the notations by identifying $T$ with the  ${\bf X}$-tree whose labeling map is simply the identity function. Moreover, if ${\bf X}' \subset {\bf X}$, we will  also identify $T |{\bf X}'$ with  the ${\bf X}'$-tree that is the restriction of $T$ (as an ${\bf X}$-tree) to ${\bf X}'$.

Now we begin to prove \lemref{pseudo}. The ``only if" part of the theorem is trivial and we will only prove the ``if" part of the statement.

 We recall that ${\bf X} = \{X_1, \dots, X_m\}$. Let $\delta$ be a pseudo-metric on ${\bf X}$ satisfying the two conditions \eqref{T4PC} and \eqref{T3PC} in display. For any four distinct points $p, q, r , s \in [m]$, given the tree structure of $T$  it must be true that $ph_T({\pi_p}, {\pi_q}) \cap ph_T({\pi_r}, {\pi_s}) = \emptyset$ for some permutation $\pi$ of $p, q, r, s$. By \eqref{T4PC}, together with the fact that $\delta$ is a pseudo-metric, $\delta$ is in fact a \emph{tree metric} (\citet[Theorem 7.2.6]{phylo}), i.e., there exists an ${\bf X}$-tree $\tilde{\mathcal{T}}= (\tilde{T}, \tilde{\phi})$ for a tree $\tilde{T} = (\tilde{V}, \tilde{E})$ and a labeling map $\tilde{\phi}: {\bf X} \rightarrow \tilde{V}$, as well as  a \emph{strictly positive} weighting function $\tilde{w}: \tilde{E} \longrightarrow \mathbb{R}_{>0}$ such that 
\[
\delta_{pq}= \left\{
  \begin{array}{lr}
    \sum_{\tilde{e} \in ph_{\tilde{T}}(\tilde{\phi}(X_p), \tilde{\phi}(X_q))} \tilde{w}(\tilde{e}) & \text{ if } \ \ \tilde{\phi}(X_p) \not= \tilde{\phi}(X_q)\\
    0 & :  \tilde{\phi}(X_p) = \tilde{\phi}(X_q)
  \end{array}
\right.
\]
for all $p, q \in [m]$.  By Theorem 6.3.5$(i)$ and Lemma 7.1.4 in \citet{phylo}, to show that $\delta$ can be induced from $T$ it suffices to show that for any ${\bf X}' \subset {\bf X}$ of size at most $4$, the two restricted ${\bf X}'$-trees $\tilde{\mathcal{T}}|{\bf X}'$ and $T|{\bf X}'$ are such  $\Sigma(\tilde{\mathcal{T}}|{\bf X}') \subset \Sigma(T|{\bf X}')$.  Note that this is trivial for $|{\bf X}'| = 1$ and $|{\bf X}'| = 2$. For $3 \leq  |{\bf X}'| \leq 4$, we  first note that 
\begin{multline} \label{3and1}
\{X_p\}| {\bf X}'\backslash \{X_p\} \in \Sigma(\tilde{\mathcal{T}}|{\bf X}') \text{ if and only if  } 
\delta_{pq} + \delta_{pr} - \delta_{qr} > 0 \text{ for all } X_q, X_r \in {\bf X}'\backslash \{X_p\}
\end{multline}
and
\begin{multline} \label{2and2}
\{X_p, X_q\} | \{X_r, X_s\}  \in  \Sigma(\tilde{\mathcal{T}}|{\bf X}')\text{ if and only if }\\
\delta_{pr} + \delta_{qs} - \delta_{pq} - \delta_{rs} > 0 \ \  (\text{ and } \delta_{ps} + \delta_{qr} - \delta_{pq} - \delta_{rs} > 0). 
\end{multline}
These characterization for the elements in $ \Sigma(\tilde{\mathcal{T}}|{\bf X}')$ can be easily checked; also see \citet[p.148]{phylo} where these characterizations are stated. To finish the proof it remains to show that, when $3 \leq  |{\bf X}'| \leq 4$, any ${\bf X}'$-split $\{X_p\}| {\bf X}'\backslash \{X_p\}$ as in \eqref{3and1} or any ${\bf X}'$-split $\{X_p, X_q\} | \{X_r, X_s\} $ as in \eqref{2and2} must also be an element of $\Sigma(T|{\bf X}')$.

First, towards a contradiction, suppose there exists an ${\bf X}'$-split $\{X_p\}| {\bf X}'\backslash \{X_p\}$ that is an element of  $\Sigma(\tilde{\mathcal{T}}|{\bf X}')$ but not an element of $\Sigma(T|{\bf X}')$. Since $\{X_p\}| {\bf X}'\backslash \{X_p\}$ is not an element of $\Sigma(T|{\bf X}')$, by considering $T|{\bf X}'$ as a tree it must be the case that the node $X_p$ has degree at least two. Then, by condition \eqref{T3PC}, there must exist two distinct $X_q$ and $X_r$ in the set ${\bf X}'\backslash \{X_p\}$ such that $\delta_{pq} + \delta_{pr} = \delta_{qr}$. But this reaches a contradiction since by \eqref{3and1} $\delta_{pq} + \delta_{pr} - \delta_{qr} > 0$ as $\{X_p\}| {\bf X}'\backslash \{X_p\}\in \Sigma(\tilde{\mathcal{T}}|{\bf X}')$. 

Similarly, suppose $\{X_p, X_q\} | \{X_r, X_s\}$ is an element of $ \Sigma(\tilde{\mathcal{T}}|{\bf X}')$ but not an element of $\Sigma(T|{\bf X}')$. Since $\{X_p, X_q\} | \{X_r, X_s\} \in  \Sigma(\tilde{\mathcal{T}}|{\bf X}')$, the two strict inequalities in  \eqref{2and2} must be true. On the other hand, if $T| {\bf X}'$ has any of the configurations in \figref{config}$(a) - (c)$,  since  $\{X_p, X_q\} | \{X_r, X_s\} \not\in  \Sigma(T|{\bf X}')$ it must be true that $ph_T( p, q) \cap ph_T(r, s) \not = \emptyset$, in which case it must lead to either $ph_T( p, r) \cap ph_T(q, s) = \emptyset$ or $ph_T( p, s) \cap ph_T(q, r) = \emptyset$, contradicting one of the  inequalities in $\eqref{2and2}$  by condition \eqref{T4PC}. If $T| {\bf X}'$ has the configuration in \figref{config}$(d)$ or $(e)$, then it must be the case that both $ph_T( p, r) \cap ph_T(q, s) $ and  $ph_T( p, s) \cap ph_T(q, r)$ are empty sets, which also contradict both inequalities in $\eqref{2and2}$  by condition \eqref{T4PC}.

\subsubsection*{Acknowledgments}

Part of this work was undertaken while Dennis Leung was a postdoc at the Chinese University of Hong Kong, and he would like to thank Professor Qi-Man Shao for some helpful discussions during that time.

\bibliographystyle{plainnat}

\bibliography{AlgTestHD_bib}

\end{document}